\documentclass{amsart}

\newtheorem{theorem}{Theorem}[section]
\newtheorem{lemma}[theorem]{Lemma}

\newtheorem{corollary}[theorem]{Corollary}

\theoremstyle{definition}

\def\R{{\mathbb R}}

\def\N{{\mathbb N}}

\def \cI {\mathcal I}

\theoremstyle{remark}

\newtheorem*{note*}{Note}
\numberwithin{equation}{section}
\makeatletter
\newcommand{\rank}{\mathop{\operator@font rank}}
\newcommand{\conv}{\mathop{\operator@font conv}}
\newcommand{\vol}{\mathop{\operator@font vol}}
\newcommand{\onetagright}{\tagsleft@false}
\makeatother
\renewcommand{\epsilon}{\varepsilon}

\usepackage{color}

\usepackage{hyperref}

\begin{document}

\title{Estimates for moments of  general measures on convex bodies}

\author{Sergey Bobkov}
\address{Department of Mathematics, University of Minnesota,
206 Church St SE, Minneapolis, MN 55455 USA}
\email{bobkov@math.umn.edu}

\author{Bo'az Klartag}
\address{Department of Mathematics, Weizmann Institute of Science,
Rehovot 76100 Israel, and
School of Mathematical Sciences, Tel Aviv University, Tel Aviv 69978}
\email{boaz.klartag@weizmann.ac.il}

\author{Alexander Koldobsky}
\address{Department of
Mathematics, University of Missouri, Columbia, MO 65211}
\email{koldobskiya@missouri.edu}
\thanks{This material is based upon work supported by the US National Science
Foundation under Grant DMS-1440140 while the authors were in residence
at the Mathematical Sciences Research Institute in Berkeley, California,
during the Fall 2017 semester.
The first- and third-named authors were supported in part by the NSF
Grants DMS-1612961 and DMS-1700036.
The second-named author was supported in part by a European Research Council
(ERC) grant.}

\subjclass{}
\keywords{}

\begin{abstract} For $p\ge 1$, $n\in \N$, and an origin-symmetric convex
body $K$ in $\R^n,$ let
$$
d_{\rm {ovr}}(K,L_p^n) =
\inf \left\{ \Big(\frac{|D|}{|K|}\Big)^{1/n}: K \subseteq D,\
D\in L_p^n \right\}
$$
be the outer volume ratio distance from $K$ to the class $L_p^n$ of the unit
balls of $n$-dimensional subspaces of $L_p.$ We prove that there exists
an absolute constant $c>0$ such that
\begin{equation}\label{dist}
\frac{c\sqrt{n}}{\sqrt{p\log\log n}}\le \sup_K d_{\rm {ovr}}(K,L_p^n)\le \sqrt{n}.
\end{equation}
This result follows from a new slicing inequality for arbitrary measures,
in the spirit of the slicing problem of Bourgain. Namely, there exists
an absolute constant $C>0$ so that for any $p\ge 1,$ any $n\in \N$,
any compact set $K \subseteq \R^n$ of positive volume, and any
Borel measurable function $f\ge 0$ on $K$,
\begin{equation}\label{slicing}
\int_K f(x)\,dx \, \le \, C \sqrt{p}\ d_{\rm ovr}(K,L_p^n)\ |K|^{1/n}
\sup_{H} \int_{K\cap H} f(x)\,dx,
\end{equation}
where the supremum is taken over all affine hyperplanes $H$ in $\R^n$.
Combining (\ref{slicing}) with a recent counterexample for the slicing problem
with arbitrary measures from \cite{KK}, we  get the lower estimate
from (\ref{dist}).

In turn, inequality (\ref{slicing}) follows from an estimate for the $p$-th
absolute moments of the function $f$
$$
\min_{\xi \in S^{n-1}} \int_K |(x,\xi)|^p f(x)\ dx \, \le \,
(Cp)^{p/2}\, d^p_{\rm {ovr}}(K,L_p^n)\ |K|^{p/n} \int_K f(x)\,dx.
$$
Finally, we prove a result of the Busemann-Petty type for these moments.
\end{abstract}

\maketitle

%%%%%%%%%%%%%%%%%%%%%%%%%%%%%%%%%%%%%%%%%%%%%%%%%%%%%%%%%%%%%%%%%%%%%%%%%%%%%%%%%%%%%%%%%%%%%%%%%%%%%%%%%%%%%%%%%%%%%%%%%%%%%%%%%%%%%%%%%%%%%%%%%%
\section{Introduction}
%%%%%%%%%%%%%%%%%%%%%%%%%%%%%%%%%%%%%%%%%%%%%%%%%%%%%%%%%%%%%%%%%%%%%%%%%%%%%%%%%%%%%%%%%%%%%%%%%%%%%%%%%%%%%%%%%%%%%%%%%%%%%%%%%%%%%%%%%%%%%%%%%%

Suppose that $K \subseteq \R^n$ ($n \geq 1$) is a centrally-symmetric convex
set of volume one (i.e., $K = -K$). Given an even continuous probability density
$f: K \rightarrow [0,\infty)$, and $p\ge 1$, can we find a direction
$\xi$ such that the $p$-th absolute moment
\begin{equation} \label{sergey-question}
M_{K,f,p}(\xi) = \int_{K} |(x,\xi)|^p\, f(x) \,dx
\end{equation}
is smaller than a constant which does not depend on $K$ and $f$?
More precisely and in a more relaxed form, let $\gamma(p,n)$ be the smallest
number $\gamma > 0$ satisfying

\begin{equation}\label{gammap}
\min_{\xi\in S^{n-1}} M_{K,f,p}(\xi) \le \gamma^p\ |K|^{p/n} \int_K f(x)\,dx
\end{equation}
for all centrally-symmetric convex bodies $K \subseteq \R^n$
and all even continuous functions $f \geq 0$ on $K$.
Here and below, we denote by $S^{n-1} = \{\xi \in \R^n: |\xi| = 1 \}$
the Euclidean unit sphere centered at the origin, and $|K|$ stands for volume
of appropriate dimension.
(Note that the continuity property of $f$ in the definition \ref{sergey-question}
is irrelevant and may easily be replaced by measurability.)
As we will see, there is a  two-sided~bound
on $\gamma(p,n)$.

\begin{theorem}\label{gamma}
With some positive absolute constants $c$ and $C$, for any $p\ge1$,
$$
\frac{c\sqrt{n}}{\sqrt{\log\log n}}\le\gamma(p,n)\le C\sqrt{pn}.
$$
\end{theorem}

To describe the way the upper bound is obtained, denote by $L_p^n$ the class
of the unit balls of $n$-dimensional subspaces of $L_p$.
Equivalently (see \cite[p. 117]{Koldobsky-book}), $L_p^n$ is the class of
all centrally-symmetric convex bodies $D$ in $\R^n$ such that there exists
a finite Borel measure $\nu_D$ on $S^{n-1}$ satisfying
\begin{equation}\label{subspLp}
\|x\|_D^p=\int_{S^{n-1}} |(x,\theta)|^p\ d\nu_D(\theta),\qquad \forall x\in \R^n.
\end{equation}
Here $\|x\|_D = \inf\{a\ge 0: x\in aD\}$ is the norm generated by $D$.
Note that $L_1^n=\Pi_n^*$ is the class of polar projection bodies which,
in particular, contains the cross-polytopes; see \cite[Ch.8]{Koldobsky-book}
for details.

For a (bounded) set $K$ in $\R^n$, define the quantity
$$
V(K,L_p^n) =
\inf \big\{|D|^{1/n}: K \subseteq D,\ D\in L_p^n \big\}.
$$
If $K$ is measurable and has positive volume, we have the relation
$$
V(K,L_p^n) = d_{\rm {ovr}}(K,L_p^n)\, |K|^{1/n},
$$
with
\begin{equation}
\label{p-ovr}
d_{\rm {ovr}}(K,L_p^n) =
\inf \left\{ \Big(\frac{|D|}{|K|}\Big)^{1/n}: K \subseteq D,\
D\in L_p^n \right\}.
\end{equation}
For convex $K$, the latter may be interpreted as the outer volume ratio
distance from $K$ to the class of unit balls of $n$-dimensional subspaces of
$L_p$. The next body-wise estimates refine the upper bound in
Theorem \ref{gamma} in terms of the $d_{\rm {ovr}}$-distance.

\begin{theorem} \label{sergey-p}
Given a probability measure $\mu$ on $\R^n$ with a compact support $K$,
for every $p\ge 1$,
$$
\min_{\xi \in S^{n-1}} \Big(\int |(x,\xi)|^p\,d\mu(x)\Big)^{1/p} \, \le \,
C\sqrt{p}\ V(K,L_p^n),
$$
where $C$ is an absolute constant. In particular, if $f$ is a non-negative
continuous function on a compact set $K \subseteq \R^n$ of positive volume, then
$$
\min_{\xi \in S^{n-1}} M_{K,f,p}(\xi) \, \le \,
(Cp)^{p/2}\, d^p_{\rm {ovr}}(K,L_p^n)\ |K|^{p/n} \int_K f(x)\,dx.
$$
\end{theorem}

In the class of centrally-symmetric convex bodies $K$ in $\R^n$,
there is a dimensional bound $d_{\rm ovr}(K,L_p^n)\le \sqrt{n}$, which
follows from John's theorem and the fact that ellipsoids belong to $L_p^n$
for all $p\ge 1$ (see \cite{John-1948} and \cite[Lemma 3.12]{Koldobsky-book}).
Hence, the second upper bound of Theorem 1.2 is more accurate in comparison
with the universal bound of Theorem 1.1.

Moreover, for several classes of centrally-symmetric convex bodies, it is
known that the distance $d_{\rm ovr}(K,L_p^n)$ is bounded by absolute constants.
These classes include duals of bodies with bounded volume ratio
(see \cite{Koldobsky-2015}) and the unit balls of normed spaces that embed
in $L_q$, $1 \le q <\infty$ (see \cite{Milman-2006, Koldobsky-Pajor}).
In the case $p=1$, they also include all unconditional convex bodies
\cite{Koldobsky-2015}.
The proofs in these papers estimate the distance from the class
of intersection bodies, but the actual bodies used there (the Euclidean ball
for $p>1$ and the cross-polytope for $p=1$) also belong to the classes
$L_p^n$, so the same arguments work for $L_p^n$.
%%%%%%%%%%%%%%

In order to prove the lower
estimate of Theorem \ref{gamma}, we first establish the connection between
question (\ref{sergey-question}) and the slicing problem for arbitrary measures.
The slicing problem of Bourgain \cite{Bou1,Bou2} asks whether
$\sup_n L_n < \infty$, where $L_n$ is the minimal positive number $L$
such that, for any centrally-symmetric convex body $K \subseteq \R^n$,
$$
|K| \, \leq \, L \max_{\xi \in S^{n-1}} |K \cap \xi^{\perp}| \,
|K|^{1/n}.
$$
Here, $\xi^\bot$ is the hyperplane in $\R^n$ passing through the origin and
perpendicular to the vector $\xi$, and we write $|K \cap \xi^{\perp}|$
for the $(n-1)$-dimensional volume.
Bourgain's slicing problem is still unsolved. The best-to-date
estimate $L_n\le C n^{1/4}$ was established by the second-named author
\cite{Klartag-2006}, removing a logarithmic term from an earlier estimate by
Bourgain \cite{Bourgain-1991}.

The slicing problem for arbitrary measures was introduced in
\cite{Koldobsky-2012} and considered in
\cite{Koldobsky-2014, Koldobsky-2015, Koldobsky-Pajor, CGL, KK}.
In analogy with the original problem, for a centrally-symmetric convex
body $K \subseteq \R^n$, let $S_{n,K}$ be the smallest positive number $S$
satisfying
\begin{equation}\label{hyperplane}
\int_K f(x)\,dx  \, \le \,
S \max_{\xi\in S^{n-1}} \int_{K\cap \xi^\bot} f(x)\,dx\ |K|^{\frac{1}{n}}
\end{equation}
for all even continuous functions $f \geq 0$ in $\R^n$
(where $dx$ on the right-hand side refers to the Lebesgue measure on
the corresponding affine subspace of $\R^n$).
It was proved in \cite{Koldobsky-2014} that
$$
S_n \, = \, \sup_{K \subseteq \R^n} S_{n,K} \, \le \, 2\sqrt{n}.
$$
However, for many classes of bodies, including intersection bodies
\cite{Koldobsky-2012} and unconditional convex bodies \cite{Koldobsky-2015},
the quantity $S_{n,K}$ turns out to be bounded by an absolute constant. In particular, if
$K$ is the unit ball of an $n$-dimensional subspace of $L_p$, $p>2$, then
$S_{n,K}\le C\sqrt{p}$ with some absolute constant $C$;
see \cite{Koldobsky-Pajor}. These results are implied
by the following estimate proved in \cite{Koldobsky-2015}:

\begin{theorem} \label{ovr} $($\cite{Koldobsky-2015}$)$
For any centrally-symmetric star body $K \subseteq \R^n$ and any even
continuous non-negative function $f$ on $K$,
$$
\int_K f(x)\,dx \, \le \, 2\ d_{\rm ovr}(K,\cI_n) \,
\max_{\xi \in S^{n-1}} \int_{K\cap \xi^\bot} f(x)\,dx \ |K|^{1/n},
$$
where $d_{\rm ovr}(K,\cI_n)$ is the outer volume ratio distance from $K$
to the class $\cI_n$ of intersection bodies in $\R^n$.
\end{theorem}

The class of intersection bodies $\cI_n$ was introduced by Lutwak \cite{Lutwak};
it can be defined as the closure in the radial metric of radial sums of
ellipsoids centered at the origin in $\R^n$.

On the other hand, it was shown in \cite{KK} that in general the constants
$S_n$ are of the order $\sqrt{n}$, up to a doubly-logarithmic term.

\begin{theorem} \label{example} $($\cite{KK}$)$
For any $n \geq 3$, there exists a centrally-symmetric convex body
$T \subseteq \R^n$ and an even, continuous probability density
$f:T \rightarrow [0, \infty)$ such that, for any affine hyperplane
$H \subseteq \R^n$,
 \begin{equation}
\int_{T \cap H} f(x)\,dx \, \leq \,
C \frac{\sqrt{\log \log n}}{\sqrt{n}} \, |T|^{-1/n},
\label{eq_649}
\end{equation}
where $C > 0$ is a universal constant.
\label{thm_639} \end{theorem}

The connection between (\ref{gammap}) and the slicing inequality
for arbitrary measures (\ref{hyperplane}) is as follows.

\begin{lemma}\label{sections-p}
Given a Borel measurable function $f \geq 0$ on $\R^n$,
for any $\xi\in S^{n-1}$ and $p>0$,
$$
2^p\,(p+1)\, \bigg(\sup_{s\in \R} \int_{(x,\xi)=s} f(x)\,dx\bigg)^p \
\int |(x,\xi)|^p\, f(x) \,dx
 \, \ge \,
\bigg(\int f(x)\,dx\bigg)^{p+1}.
$$
\end{lemma}

\vskip5mm
If $f$ is defined on a set $K$ in $\R^n$, we then have
$$
2^p\,(p+1)\, \bigg(\sup_{s\in \R} \int_{K\cap \{(x,\xi)=s\}} f(x)\,dx\bigg)^p \
M_{K,f,p}(\xi)
 \, \ge \,
\bigg(\int_K f(x)\,dx\bigg)^{p+1}.
$$
The lower bound in Theorem \ref{gamma} thus follows, by combining
the above inequality with (\ref{gammap}) and Theorem \ref{example}.

\begin{corollary} With some positive absolute constants $c$ and $C$,
for every $p\ge 1$,
$$
\frac{c\sqrt{n}}{\sqrt{\log\log n}}\le S_n\le C\gamma(p,n).
$$
\end{corollary}

Lemma \ref{sections-p}, in conjunction with Theorem \ref{sergey-p}, leads
to a new slicing inequality. In the case of volume, where $f\equiv 1$,
this inequality was established earlier by Ball \cite{Ball-1991} for
$p=1$ and by Milman \cite{Milman-2006} for arbitrary $p$.

\begin{theorem}\label{slicing-p}
Let $f \geq 0$ be a Borel measurable function on a compact
set $K \subseteq \R^n$ of positive volume. Then, for any $p>2$,
$$
\int_K f(x)\,dx \, \le \, C \sqrt{p}\ d_{\rm ovr}(K,L_p^n)\ |K|^{1/n}
\sup_{H} \int_{K\cap H} f(x)\,dx,
$$
where the supremum is taken over all affine hyperplanes $H$ in $\R^n$,
and $C$ is an absolute constant.
\end{theorem}

Theorem \ref{slicing-p} also holds for $1\le p\le 2$, but in this case it
is weaker than Theorem \ref{ovr}, because the unit ball of every finite
dimensional subspace of $L_p,\ 0<p\le 2,$ is an intersection body;
see \cite{Koldobsky-1998}. However, for $p>2$ the unit balls of  subspaces
of $L_p$ are not necessarily intersection bodies. For example the unit balls
of $\ell_p^n$ are not intersection bodies if $p>2$, $n\ge 5$;
see \cite[Th. 4.13]{Koldobsky-book}. So the result of Theorem \ref{slicing-p}
is new for $p>2$, and generalizes the estimate from \cite{Koldobsky-Pajor}
in the case where $K$ itself belongs to the class $L_p^n$.

Theorem \ref{slicing-p} gives another reason to estimate
the outer volume ratio distance $d_{\rm ovr}(K,L_p^n)$
from an arbitrary symmetric convex body to the class of unit balls of subspaces
of $L_p$. As mentioned before,
$$d_{\rm ovr}(K,L_p^n)\le \sqrt{n},
$$
uniformly over all centrally-symmetric convex bodies $K$ in $\R^n$.
Surprisingly, the corresponding lower estimates seem to be missing in
the literature. Combining Theorems \ref{slicing-p} and  \ref{example}, we get
a lower estimate which shows that $\sqrt{n}$ is optimal up to a doubly-logarithmic term
with respect to the dimension $n$ and a term depending on the power $p$ only.

\begin{corollary} \label{distance-p}
There exists a centrally-symmetric convex body $T \subseteq \R^n$ such that
$$
d_{\rm ovr}(T, L_p^n) \, \geq \, c\, \frac{\sqrt{n}}{\sqrt{p\log \log n}}
$$
for every $p\ge 1$, where $c > 0$ is a universal constant.
\end{corollary}

We end the Introduction  with a comparison result for the
quantities $M_{K,f,p}(\xi)$. For $p\ge 1$, introduce the Banach-Mazur distance
$$
d_{BM}(M,L_p^n) =
\inf\left\{a \geq 1:\, \exists\, D\in L_p^n \ {\rm such \ that} \
D\subset M\subset aD \right\}
$$
from a star body $M$ in $\R^n$ to the class $L_p^n$.
Recall that $L_p^n$ is invariant with respect to linear transformations.
By John's theorem, if $M$ is origin-symmetric and convex, then
$d_{BM}(M,L_p^n)\le \sqrt{n}$. We prove the following:

\begin{theorem}\label{bp}
Let $K$ and $M$ be origin-symmetric star bodies in $\R^n$, and let $f \geq 0$
be an even continuous function on $\R^n$.
Given $p\ge 1$, suppose that for every $\xi\in S^{n-1}$
\begin{equation}\label{cond1}
\int_K |(x,\xi)|^p\, f(x)\,dx \, \le \, \int_M |(x,\xi)|^p\, f(x)\, dx.
\end{equation}
Then
$$
\int_K f(x)\,dx \, \le \, d^p_{BM}(M,L_p^n) \int_M f(x)\,dx.
$$
\end{theorem}

\vskip2mm
This result is in the spirit  of the isomorphic Busemann-Petty problem for
arbitrary measures proved in \cite{KZ}: with the same notations, if
$$
\int_{K\cap \xi^\bot} f(x)\,dx \, \le \,
\int_{M\cap \xi^\bot} f(x)\,dx,\qquad \forall \xi\in S^{n-1},
$$
then
$$
\int_K f(x)\,dx \, \le \, d_{BM}(K,\cI_n) \int_M f(x)\,dx.
$$
We refer the reader to \cite[Ch.5]{Koldobsky-book} for more about
the Busemann-Petty problem.

Throughout this paper, we write  $a\sim b$ when $ca\le b \le Ca$ for
some absolute constants $c,C.$ A convex body $K$ in $\R^n$
is a compact, convex set with a non-empty interior.
The standard scalar product between $x,y \in \R^n$ is denoted by $(x,y)$
and the Euclidean norm of $x \in \R^n$ by $|x|$.
We write $\log$ for the natural logarithm.

\section{Proofs}
In this section we prove Theorem \ref{sergey-p}, Lemma \ref{sections-p} and
Theorem \ref{bp}. The other results of this paper will follow as explained
in the Introduction.

Given a compact set $K \subseteq \R^n$ and $x \in \R^n$, put
$$
\|x\|_K = \min\{a\ge 0:\ x\in aK\},
$$
if $x\in aK$ for some $a \geq 0$, and $\|x\|_K = \infty$
in the other case. For star bodies, it represents the usual
{\it Minkowski functional} associated with $K$.

\medbreak
\noindent{\bf Proof of Theorem \ref{sergey-p}.}
Let $D \subseteq \R^n$ be the unit ball of
an $n$-dimensional subspace of $L_p$, so that the relation (1.3) holds
for some measure $\nu_D$ on the unit sphere $S^{n-1}$.
Then, integrating the inequality
$$
\min_{\theta\in S^{n-1}} \int_K |(x,\theta)|^p\ d\mu(x) \le
\int_K |(x,\xi)|^p\ d\mu(x) \qquad (\xi \in S^{n-1})
$$
over the variable $\xi$ with respect to $\nu_D$, we get the relation
$$
\nu_D(S^{n-1})\min_{\theta\in S^{n-1}} \int_K |(x,\theta)|^p\ d\mu(x) \, \le \,
\int_K \|x\|_D^p\,d\mu(x).
$$
In the case $K\subseteq D$, we have $\|x\|_D \le \|x\|_K \leq 1$ on $K$,
so that the last integral does not exceed $\mu(K) = 1$, and thus
\begin{equation}
\nu_D(S^{n-1})\min_{\theta\in S^{n-1}} \int_K |(x,\theta)|^p\ d\mu(x)
 \, \le \, 1.
\end{equation}

In order to estimate the left-hand side of (2.1) from below,
we represent the value $\nu_D(S^{n-1})$ as the integral
$\int_{S^{n-1}} |x|^p\, d\nu_D(x)$ and apply
the well-known formula
$$
|x|^p \, = \,
\frac{\Gamma(\frac{p+n}2)}{2\pi^{\frac{n-1}2}\,\Gamma(\frac{p+1}2)}\,
\int_{S^{n-1}} |(x,\theta)|^p\,d\theta, \qquad x\in \R^n
$$
(see for example \cite[Lemma 3.12]{Koldobsky-book}). Using (\ref{subspLp}),
this yields the representation
\begin{eqnarray*}
\nu_D(S^{n-1})
 & = &
\frac{\Gamma(\frac{p+n}2)}{2\pi^{\frac{n-1}2}\Gamma(\frac{p+1}2)}\,
\int_{S^{n-1}} \int_{S^{n-1}} |(x,\theta)|^p\ d\theta\ d\nu_D(x) \\
 & = &
\frac{\Gamma(\frac{p+n}2)}{2\pi^{\frac{n-1}2}\Gamma(\frac{p+1}2)}\,
\int_{S^{n-1}} \|\theta\|_D^p\, d\theta.
\end{eqnarray*}
The last integral may be related to the volume of $D$, by using
the polar formula for the volume of $D$,
$$
n\,|D| = \int_{S^{n-1}} \|\theta\|_D^{-n}\,d\theta
 = s_{n-1} \int_{S^{n-1}} \|\theta\|_D^{-n}\,d\sigma_{n-1}(\theta),
$$
where $\sigma_{n-1}$ denotes the normalized Lebesgue measure on $S^{n-1}$
and $s_{n-1} = \frac{2\pi^{\frac{n}{2}}}{\Gamma(\frac{n}{2})}$ is its
$(n-1)$-dimensional volume. Namely, by Jensen's inequality, we have
$$
\int \|\theta\|_D^{-n}\,d\sigma_{n-1}(\theta) \, \geq \,
\bigg(\int \|\theta\|_D^p\,d\sigma_{n-1}(\theta)\bigg)^{-\frac{n}{p}},
$$
or equivalently
$$
\int \|\theta\|_D^p\,d\theta \, \geq \,
s_{n-1}^{\frac{p+n}{n}}\,(n\,|D|)^{-\frac{p}{n}}.
$$
Thus,
\begin{eqnarray*}
\nu_D(S^{n-1})
 & \ge &
\frac{\Gamma(\frac{p+n}2)\,s_{n-1}^{\frac{p+n}{n}}}{2\pi^{\frac{n-1}2}\,
\Gamma(\frac{p+1}2)\,n^{\frac{p}{n}}\,|D|^{\frac{p}{n}}} \\
 & = &
\sqrt{\pi}\,\frac{\Gamma(\frac{p+n}{2})}{\Gamma(\frac{p+1}{2})\,
\Gamma(\frac{n}{2})}\,
\Big(\frac{s_{n-1}}{n\,|D|}\Big)^{\frac{p}{n}}
 \ \ge \
\frac{c^p}{\Gamma(\frac{p+1}{2})\,|D|^{\frac{p}{n}}},
\end{eqnarray*}
where $c>0$ is an absolute constant. Here we used the well-known asymptotic
relation $\sqrt{n}\,s_{n-1}^{\frac{1}{n}} \rightarrow c_0$ as
$n \rightarrow \infty$, for some absolute $c_0>0$, as well as the 
estimate $\Gamma(\frac{p+n}{2}) / \Gamma( \frac{n}{2}) \geq ( cn)^{p/2}$.

Applying this lower estimate on the left-hand side of (2.1), we get
$$
\min_{\theta\in S^{n-1}} \int_K |(x,\theta)|^p\ d\mu(x) \, \leq \,
C^p\, \Gamma\Big(\frac{p+1}{2}\Big)\,|D|^{\frac{p}{n}}.
$$
It remains to take the minimum over all admissible $D$ and note that
$\Gamma\left(\frac{p+1}2\right)^{1/p} \leq c\sqrt{p}$ for $p \ge 1$.
\qed
\medbreak

To prove Lemma \ref{sections-p}, we need the following simple assertion.

\begin{lemma}\label{p>q}
Given a measurable function $g:\R \rightarrow [0,1]$, the function
$$
q\mapsto \bigg(\frac{q+1}{2} \int_{-\infty}^\infty
|t|^q\,g(t)\, dt\bigg)^{\frac{1}{q+1}}
$$
is non-decreasing on $(-1,\infty)$.
\end{lemma}

\noindent{\bf Proof.}
The standard argument is similar to the one used in the proof of Lemma 2.4 in
\cite{Klartag-2005}. Given $-1<q<p$, let $A>0$ be defined by
$$
\int_{-\infty}^\infty |t|^q\, g(t)\, dt = \int_{-A}^A |t|^q\, dt =
\frac{2}{q+1}\,A^{q+1}.
$$
Using
$$
|t|^p \leq A^{p-q}\, |t|^q \ \ (|t| \leq A) \quad {\rm and} \quad
|t|^p \geq A^{p-q}\, |t|^q \ \ (|t| \geq A),
$$
together with the assumption $0 \leq g \leq 1$, we then have
\begin{eqnarray*}
\hskip5mm
\int_{|t| \leq A} (1-g(t))\,|t|^p\, dt - \int_{|t| > A} g(t)\,|t|^p\, dt
 & \\
 &
\hskip-45mm \le \
A^{p-q}\,\Big(\int_{|t| \leq A} (1-g(t))\,|t|^q\,dt -
\int_{|t| > A} g(t)\,|t|^q\,dt\Big)
 \, = \, 0.
\end{eqnarray*}
Hence
$$
\int_{-\infty}^\infty g(t)\,|t|^p\, dt \, \ge \,
\int_{-A}^A |t|^p\, dt \, = \, \frac{2}{p+1}\,A^{p+1},
$$
that is,
$$
\bigg(\frac{p+1}{2}
\int_{-\infty}^\infty g(t)\,|t|^p\, dt\bigg)^{\frac{1}{p+1}}
\, \ge \, A \, = \,
\bigg(\frac{q+1}{2}
\int_{-\infty}^\infty g(t)\,|t|^q\, dt\bigg)^{\frac{1}{q+1}}.
$$
\qed

\medbreak
\noindent{\bf Proof of Lemma \ref{sections-p}.}
One may assume that $f$ is integrable.
For $t \in \R$, introduce the hyperplanes $H_t = \{(x,\xi)=t\}$.
Since $f$ is Borel measurable on $\R^n$, the function
$$
g(t) =
\frac{\int_{H_t} f(x)\,dx}{\sup_s \int_{H_s} f(x)\,dx}
$$
is Borel measurable on the line and
satisfies $\|g\|_\infty = 1$. By Fubini's theorem,
$$
\int_{-\infty}^\infty |t|^p\,g(t)\,dt
 \, = \,
\frac{\int |(x,\xi)|^p\,f(x)\,dx}{\sup_s \int_{H_s} f(x)\,dx},
$$
$$
\int_{-\infty}^\infty g(t)\,dt \, = \,
\frac{\int f(x)\,dx}{\sup_s \int_{H_s} f(x)\,dx}.
$$
Applying Lemma \ref{p>q} to the function $g$ with $q=0$ and $p$,
we get
$$
\frac{1}{2} \int_{-\infty}^\infty g(t)\, dt \, \leq \,
\bigg(\frac{p+1}{2} \int_{-\infty}^\infty
|t|^p\,g(t)\, dt\bigg)^{\frac{1}{p+1}},
$$
which in our case becomes
$$
\Big(\int f(x)\,dx\Big)^{p+1} \,\leq \, (p+1)\,
\Big(2\,\sup_s \int_{H_s} f(x)\,dx\Big)^p \int |(x,\xi)|^p\,f(x)\,dx.
$$
\qed

\smallbreak
\noindent{\bf Proof of Theorem \ref{bp}.}
Let $D\in L_p^n$ be such that the distance $d_{\rm {ovr}}(M,L_p^n)$ is
almost realized, i.e., for small $\delta>0$, suppose that
$D \subseteq M \subseteq (1+\delta)\,d_{BM}(M,L_p^n)\,D$.

Integrating both sides of (\ref{cond1}) over $\xi \in S^{n-1}$
with respect to the measure $\nu_D$ from (\ref{subspLp}),
we get
$$
\int_K \|x\|_D^p\, f(x)\,dx \, \le \, \int_M \|x\|_D^p\, f(x)\, dx.
$$
Equivalently, using the integrals in spherical coordinates, we have
$$
0 \, \le \,
\int_{S^{n-1}} \|\theta\|_D^p\,
\bigg(\int_{\|\theta\|_K^{-1}}^{\|\theta\|_M^{-1}} r^{n+p-1}
f(r\theta)\,dr\bigg) \, d\theta
 \, = \,
\int_{S^{n-1}} \frac{\|\theta\|_D^p}{\|\theta\|_M^p}\ I(\theta) \,d\theta,
$$
where
$$
I(\theta) \, = \, \|\theta\|_M^p
\int_{\|\theta\|_K^{-1}}^{\|\theta\|_M^{-1}} r^{n+p-1} f(r\theta)\,dr.
$$
For $\theta \in S^{n-1}$ such that $\|\theta\|_K \geq \|\theta\|_M$, the latter quantity is
non-negative, and one may proceed by writing
\begin{eqnarray*}
I(\theta)
 & = &
\int_{\|\theta\|_K^{-1}}^{\|\theta\|_M^{-1}}
\Big(\|\theta\|_M^p - r^{-p}\Big)\, r^{n+p-1} f(r\theta)\, dr
+
\int_{\|\theta\|_K^{-1}}^{\|\theta\|_M^{-1}} r^{n-1} f(r\theta)\, dr \\
 & \le &
\int_{\|\theta\|_K^{-1}}^{\|\theta\|_M^{-1}} r^{n-1} f(r\theta)\, dr.
\end{eqnarray*}
But, in the case $\|\theta\|_K \leq \|\theta\|_M$, we have
$$
-I(\theta) \, = \,
\|\theta\|_M^p \int_{\|\theta\|_M^{-1}}^{\|\theta\|_K^{-1}} r^p\,
r^{n-1} f(r\theta)\,dr \, \geq \,
\int_{\|\theta\|_M^{-1}}^{\|\theta\|_K^{-1}} r^{n-1} f(r\theta)\,dr,
$$
which is the same upper bound on $I(\theta)$ as before.
Thus,
$$
0 \, \le \,
\int_{S^{n-1}} \frac{\|\theta\|_D^p}{\|\theta\|_M^p}\
\bigg(\int_{\|\theta\|_K^{-1}}^{\|\theta\|_M^{-1}} r^{n-1}
f(r\theta)\,dr\bigg) \, d\theta,
$$
that is,
$$
\int_{S^{n-1}} \frac{\|\theta\|_D^p}{\|\theta\|_M^p}\
\bigg(\int_0^{\|\theta\|_K^{-1}} r^{n-1} f(r\theta)\,dr\bigg) \, d\theta
\le
\int_{S^{n-1}} \frac{\|\theta\|_D^p}{\|\theta\|_M^p}\
\bigg(\int_0^{\|\theta\|_M^{-1}} r^{n-1} f(r\theta)\,dr\bigg) \, d\theta.
$$

Now, by the choice of $D$,
$$
\|\theta\|_M \, \le \, \|\theta\|_D \, \le \,
(1+\delta)\,d_{BM}(M,L_p^n)\|\theta\|_M
$$
for every $\theta\in S^{n-1}$. Hence
\begin{eqnarray*}
\int_K f(x)\,dx
 & = &
\int_{S^{n-1}}
\bigg(\int_{0}^{\|\theta\|_K^{-1}} r^{n-1} f(r\theta)\,dr\bigg)\,d\theta \\
 & \le &
\int_{S^{n-1}} \frac{\|\theta\|_D^p}{\|\theta\|_M^p}\
\bigg( \int_{0}^{\|\theta\|_K^{-1}} r^{n-1} f(r\theta)\,dr\bigg)\,d\theta \\
 & \le &
\int_{S^{n-1}} \frac{\|\theta\|_D^p}{\|\theta\|_M^p}\
\bigg(\int_{0}^{\|\theta\|_M^{-1}} r^{n-1} f(r\theta)\,dr\bigg)\,d\theta \\
 & \le &
(1+\delta)\,d^p_{BM}(M,L_p^n) \int_{S^{n-1}}
\bigg(\int_{0}^{\|\theta\|_M^{-1}} r^{n-1} f(r\theta)\,dr\bigg)\,d\theta \\
 & = &
(1+\delta)\,d^p_{BM}(M,L_p^n) \int_M f(x)\,dx.
\end{eqnarray*}
Sending $\delta$ to zero, we get the result. \qed

{
}

\begin{thebibliography}{99}
\addcontentsline{toc}{section}{References}
\setlength{\itemsep}{1pt}


\bibitem{Ball-1991} K.~Ball, \textit{Normed spaces with a weak Gordon-Lewis property},
Lecture Notes in Math. {\bf 1470}, Springer, Berlin (1991), 36--47.

\bibitem{Bou1} J.~Bourgain, {\it On high-dimensional maximal functions associated to convex bodies}, Amer. J. Math., {108}, (1986), 1467–-1476.

\bibitem{Bou2} J.~Bourgain, {\it  Geometry of Banach spaces and harmonic analysis },
Proceedings of the International Congress of Mathematicians, Vol. 1, 2 (Berkeley, Calif., 1986), Amer. Math. Soc., Providence, RI, (1987), 871–-878.

\bibitem{Bourgain-1991} {J.~Bourgain}, {\it On the distribution of polynomials
on high-dimensional convex sets},
Geom. aspects of Funct. Anal. (GAFA seminar notes), Israel Seminar, Springer Lect.
Notes in Math. { 1469} (1991),
127--137.

\bibitem{CGL} {G.~Chasapis, A.~Giannopoulos and D-M.~Liakopoulos},
{\it Estimates for measures of lower dimensional sections of convex bodies}, Adv. Math. { 306} (2017), 880--904.

\bibitem{John-1948} {F.~John}, \textit{Extremum problems with inequalities as subsidiary conditions}, Courant
Anniversary Volume, Interscience, New York (1948), 187-204.

\bibitem{Klartag-2005} B.~Klartag, {\it An isomorphic version of the slicing problem},
J. Funct.  Anal.  218 (2005), 372 -- 394.

\bibitem{Klartag-2006}  {B.~Klartag}, {\it On convex perturbations with a bounded
isotropic constant}, Geom. Funct. Anal. (GAFA) { 16}  (2006), 1274--1290.

\bibitem{KK} B.~Klartag and A.~Koldobsky, {\it An example related to the slicing inequality for general measures},
J. Funct.  Anal., in press.

\bibitem{Koldobsky-1998} {A.~Koldobsky}, {\it Intersection bodies, positive definite
distributions and the Busemann-Petty problem}, Amer. J. Math. { 120} (1998),
827--840.

\bibitem{Koldobsky-book} A.~Koldobsky, {\it Fourier analysis in convex geometry},
Amer. Math. Soc., Providence RI, 2005.

\bibitem{Koldobsky-2012} {A.~Koldobsky}, {\it A hyperplane inequality for measures of convex bodies in
$\R^n, n\le 4$}, Discrete Comput. Geom. { 47} (2012), 538--547.

\bibitem{Koldobsky-2014} {A.~Koldobsky}, {\it A $\sqrt{n}$ estimate for measures of hyperplane sections of convex bodies,} Adv. Math. { 254} (2014), 33--40.

\bibitem{Koldobsky-2015} A.~Koldobsky, {\it Slicing inequalities for measures of convex bodies,}
Adv. Math. 283 (2015), 473--488.

\bibitem{Koldobsky-Pajor} A.~Koldobsky and A.~Pajor, {\it A remark on measures of sections of $L_p$-balls},
Geom. aspects of Funct. Anal. (GAFA seminar notes), Israel Seminar, Springer Lect. Notes in Math. 2169 (2017), 213--220.

\bibitem{KZ} A.~Koldobsky and A.~Zvavitch, {\it An isomorphic version of the Busemann-Petty problem for arbitrary measures}, Geom. Dedicata 174 (2015), 261–-277.

\bibitem{Lutwak} { E.~Lutwak}, {\it Intersection bodies and dual mixed volumes},
Adv. Math. {71} (1988), 232--261.

\bibitem{Milman-2006} E.~Milman, \textit{Dual mixed volumes and the slicing problem},
Adv. Math. {\bf 207} (2006), 566--598.

\end{thebibliography}
\end{document}